\documentclass[12pt,a4paper]{article}
\usepackage{amsthm,amssymb,amsmath,mathrsfs,graphicx}
\usepackage[pagebackref=false,colorlinks,linkcolor=black,citecolor=magenta]{hyperref}
\allowdisplaybreaks[1]

\newtheorem{theorem}{\bf Theorem}

\newtheorem{lemma}{\bf Lemma}[section]
\newtheorem{example}[lemma]{\bf Example}
\newtheorem{corollary}[lemma]{\bf Corollary}
\newtheorem{remark}{\bf Remark}

\renewcommand{\proof}{\noindent{\it\textbf{Proof.}}\ \ }

\newcommand{\Cay}{{\rm Cay}}
\newcommand{\SC}{{\rm SC}}

\newcommand{\eqd}{$\hfill \blacksquare$}

\newcommand{\Spec}{{\rm Spec}}
\newcommand{\Irr}{{\rm Irr}}

\date{}

\title {
Perfect State Transfer on semi-Cayley graphs over abelian groups
}

\author{Majid Arezoomand $^{a}$ \footnote{E-mail:
arezoomand@lar.ac.ir, arezoomandmajid@gmail.com}\\
{\small\em$^a$ University of Larestan, Lar, 74317-16137, Iran}
    }

\begin{document}
\maketitle
\begin{abstract}
In this paper, we consider the problem on the existence of perfect state transfer (PST for short) on semi-Cayley graphs over abelian groups (which are not necessarily regular), i.e on the graphs having semiregular and abelian subgroups of automorphisms with two orbits of equal size.  We stablish a characterization of semi-Cayley graphs over abelian groups having PST. As a result, we give a characterization of Cayley graphs over groups with an abelian subgroup of index $2$ having PST, which improves the earlier results on Cayley graphs over abelian groups, dihedral groups and dicyclic group and determines Cayley graphs over generalized dihedral groups and generalized dicyclic groups having PST.
\end{abstract}

\medskip \noindent
Keywords: Perfect state transfer, Semi-Cayley graph, Eigenvalue of a graph\\
Mathematics Subject Classification: 05C25, 81P45, 15A18.

\section{Introduction}

Let $\Gamma$ be a simple undirected graph with adjacency matrix $A$. The continuous-time quantum walk on $\Gamma$ is defined through the time-dependent unitary matrix
\begin{eqnarray*}
H(t)=\exp(-{\bf i}tA)=\sum_{k=0}^\infty \frac{(-{\bf i})^kt^kA^k}{k!}, ~~~0<t\in\Bbb R,~{\bf i}=\sqrt{-1}.
\end{eqnarray*}
This concept first introduced by Farhi and Gutman by a motivation of the Schr\"{o}dinger equation. They used it as
a paradigm to design efficient quantum algorithms \cite{FG}.

The significance of the study of quantum state transfer lies in its applications to the theory
of Quantum Information and Computation. In fact, one of the cruical ingredients in most of quantum information
 processing protocols is the transfer of a quantum state from one location to another location. Quantum spin network is an example of
physical systems that can serve as a quantum channel.  By considering the networks as graphs, in algebraic graph theory, one of the main questions is to find a characterization of graphs having perfect state transfer.  We say that a graph $\Gamma$ has a {\it perfect state transfer} (PST for short) from the vertex $u$ to the vertex $v$ at time $t$ if the
$(u,v)$-entry of $H(t)$, denoted by $H_{u,v}(t)$, has absolute value $1$. If $|H_{u,u}(t)|=1$ then we say that $\Gamma$ is periodic
at $u$ with period $t$.  $\Gamma$ is called periodic if it is periodic at all vertices with period $t$.

Recently, the problem of characterization of graphs having PST has been one of the interesting topics in algebraic graph theory.  For a survey on this topic and its related question, up to 2011, we refer the reader to \cite{G1,G2,G3}.
Godsil has showed that if perfect state transfer
occurs in a graph, then the square of its spectral radius is either an integer or lies in a quadratic
extension of the rationales \cite{G3}. Moreover, existence of PST has been obtained for some
particular classes of graphs, for example, distance-regular graphs \cite{CGG,C}, 
Hadamard diagonalizable graphs \cite{J} and complete bipartite
graphs \cite{SS}. Also some authors studied the PST problem of Cayley graphs due to their nice algebraic structures.  A characterization of Cayley graphs over abelian group, dihedral groups with conjugate-closed Cayley set, dihedral groups with non-conjugate-closed Cayley set, dicyclic groups and also over semidihedral groups are given in \cite{TFC, CF, CCL,ASGh, LCWW}, respectively. 
Despite many results on PST on graphs, there are no known constructions of PST on
non-abelian Cayley graphs in literature. 

The purpose of this paper is two-fold: First we present an explicit  general method to find the eigenvalues and eigenvectors of the semi-Cayley graphs over abelian groups which facilitates  the computations of PST. Many families of well-known graphs are examples of semi-Cayley graphs over abelian groups, such as bicirculants (including I-graphs, the generalized Petersen graphs, rose window graphs, Taba\v{c}jn graphs, $n$-gonal $n$-cone graphs), $n$-sunlet graphs, and the join graph of two Cayley graphs over the same group, also all Cayley graphs over groups having a subgroup of index $2$. So we provide a way to characterize
a large family of graphs having PST.
Secondly, a vital
step for making further progress on existence of PST on Cayley graphs over non-abelian groups. Since any Cayley graph over a group
$G$ having  subgroup $H$ of index $2$ is a semi-Cayley graph over $H$, the main results of \cite{TFC, CF, CCL,ASGh, LCWW} are direct consequences of our results.
Also as a result, we investigate the existence of PST on Cayley graphs over generalized dihedral groups and generalized dicyclic groups. 

The rest of the current work is organized as follows. In  Section \ref{sec:pre}, we present some comments and notations used in the next sections. In Section \ref{sec:spectra}, we give exact formula for the eigenvalues and eigenvectors of semi-Cayley graphs over abelian groups. Finally,the existence of PST on semi-Cayley graphs over abelian groups is explored in Section \ref{sec:pst} and as a result we improve some earlier results. Also, we characterize all Cayley graphs over generalized dihedral groups and generalized dicyclic groups having PST.

\section{Preliminaries and notations}\label{sec:pre}
Here we recall some of the standard facts and notation used in this paper. Our notation for 
representations and characters of finite groups is due to the notation introduced in \cite{JL}. Throughout this paper, graphs are finite, simple and undirected, and also groups are
finite. 
\subsection{Spectrum of a graph and PST}
Let $\Gamma$ be an undirected graph with adjacency matrix $A$ and $H(t)=\exp(-{\bf i}tA)$. The eigenvalues of $A$ are called {\it (adjacency) eigenvalues} of $\Gamma$. $\Gamma$ is called {\it integral} if all of its eigenvalues are integers.
Let $\Gamma$ has $m$ vertices and $\lambda_1\geq\lambda_2\geq\ldots\geq \lambda_m$ be the eigenvalues of $\Gamma$, $v_i$
is a column eigenvector of $\lambda_i$ and $\{v_1,\ldots,v_m\}$ is an orthonormal basis of $\Bbb C^m$. Then, since $A$ is
a real symmetric matrix, by the Spectral Decomposition Theorem  \cite[Theorem 3.4]{matrix}, 
\[H(t)=\sum_{i=1}^m\exp(-{\bf i}t\lambda_i)E_i,\]
 where $E_i=v_i.v_i^*$, where $v_i^*$ is the conjugate transpose of $v_i$. So $\Gamma$ has a PST from the vertex $u$ to the vertex $v$ if and only if the the absolute value of the $(u,v)$-entry of $\sum_{i=1}^m\exp(-{\bf i}t\lambda_i)E_i$ is equal to $1$.

\subsection{Cayley graphs and semi-Cayley graphs}
Given a group $G$ and an inverse-closed subset of $G$, the Cayley graph of $G$ with respect to the Cayley set $S$, denoted by $\Cay(G,S)$, is a graph
with vertex set $G$ and edge-set $\{\{g,sg\}\mid s\in S\}$. Clearly $\Cay(G,S)$ is a vertex-transitive  $|S|$-regular graph.

As a generalization of Cayley graphs, see \cite[page 2]{AT}, semi-Cayley graphs are introduced. Let $G$ be a finite group, $R, L$ and $S$ be subsets of $G$ such that $R$ and $L$ are inverse-closed subsets not containing the identity element of $G$. The {\it semi-Cayley graph} over $G$ with respect to $R$, $L$ and $S$, denoted by $\SC(G,R,L,S)$ is
an undirected graph with vertex set $\{(g,0),(g,1)\mid g\in G\}$ and edge set consists of three sets:
 \begin{eqnarray*}
&&\{\{(x,0),(y,0)\}\mid yx^{-1}\in R\}\ \ \  (\textrm{right edges}),\\
&&\{\{(x,1),(y,1)\}\mid yx^{-1}\in L\}\ \ \   (\textrm{left edges}),\\
&&\{\{(x,0),(y,1)\}\mid yx^{-1}\in S\} \ \ \  (\textrm{spoke edges}).
\end{eqnarray*}
Clearly $\SC(G,R,L,S)$ is a regular graph if and only if $|R|=|L|$. By the following result, every Cayley graph over a group having a subgroup of index $2$ is a semi-Cayley graph. Note that the converse is not true, since as an example the Petersen graph is a semi-Cayley graph over the cyclic group $\Bbb Z_5$ and it not a Cayley graph. We will use this fact frequently without referring it.

\begin{lemma}{\rm (see \cite[Lemma 8]{AT})} Let $\Gamma=\Cay(G,T)$, and $H$ be a subgroup of $G$ of index $2$ with left transversal $\{1,x\}$. Then $\Gamma\cong \SC(H,R,L,S)$, where $R=H\cap T$, $L=H\cap xTx^{-1}$ and $S=H\cap xT$.
\end{lemma}  
\subsection{Irreducible characters of finite abelian groups}
In this paper, we need irreducible characters of finite abelian groups. For the definitions and basic facts about representations and characters of finite groups, not defined here, we refer the reader to \cite{JL}. Let $G$ be an abelian
group of order $n$. It is well-known that every finite abelian group is a direct product of cyclic groups. Let $\Bbb Z_r$ be 
a cyclic group of order $r$. Then we may assume that $G=\langle g_1\rangle\times\langle g_2\rangle\times\ldots\times\langle g_k\rangle$
where $k\geq 1$ and $n=n_1n_2\ldots n_k$ and $\langle g_l\rangle\cong\Bbb Z_{n_l}$. Let $\Irr(G)$ be the set of all inequivalent complex irreducible characters of $G$. Then, by \cite[Theorem 9.8]{JL},
$\Irr(G)=\{\chi_{j_1,\ldots,j_k}\mid 0\leq j_l\leq n_{l}, l=1,\ldots,k\}$, where 
\[\chi_{j_1,\ldots,j_k}:G\rightarrow\Bbb C,~~~(g_1^{i_1},\ldots,g_k^{i_k})\mapsto \exp(\frac{2\pi {\bf i}j_1i_1}{n_1})\ldots\exp(\frac{2\pi {\bf i}j_ki_k}{n_k}).\]
The character $\chi_{0,\ldots,0}$ is the trivial character and in the rest of paper, we denote it by $\chi_1$. It is clear
that for every $g\in G$ and $\chi\in\Irr(G)$, we have $|\chi(g)|=1$. For a character of $G$ and a subset $X$ of $G$, we denote
$\sum_{x\in X}\chi(x)$ with $\chi(X)$ and if $X=\emptyset$ then we put $\chi(X)=0$.

\section{Spectra of semi-Cayley graphs over abelian groups}\label{sec:spectra}

Let $\Gamma=\SC(G,R,L,S)$ be a semi-Cayley graph over an abelian group $G$. Let $G=\{g_1,\ldots,g_n\}$ and
$\Irr(G)=\{\chi_1=1,\ldots,\chi_n\}$
be the set of all inequivalent irreducible characters of $G$. For $1\leq i\leq n$, we define the matrix
\begin{eqnarray*}
A_i=\begin{bmatrix}
\chi_i(R)&\chi_i(S)\\
\overline{\chi_i(S)}&\chi_i(L)
\end{bmatrix},
\end{eqnarray*}
where $\overline{\chi_i(S)}$ is the complex conjugate of $\chi_i(S)$. Note that since $R$ and $L$ are inverse-closed, $\chi_i(R)$ and $\chi_i(L)$ are real numbers by \cite[Proposition 13.9]{JL}. Then eigenvalues of $A_i$ are
\begin{eqnarray*}
&&\lambda_i^+=\frac{1}{2}\big(\chi_i(R)+\chi_i(L)+\sqrt{(\chi_i(R)-\chi_i(L))^2+4|\chi_i(S)|^2}\big),
\end{eqnarray*}
and
\begin{eqnarray*}
&&\lambda_i^-=\frac{1}{2}\big(\chi_i(R)+\chi_i(L)-\sqrt{(\chi_i(R)-\chi_i(L))^2+4|\chi_i(S)|^2}\big).
\end{eqnarray*}

Let $u_i^+=(a_i^+,b_i^+)$ and  $u_i^-=(a_i^-,b_i^-)$ be an eigenvector of $A_i$ corresponding to $\lambda_i^+$ and $\lambda_i^-$, respectively. Then

\begin{eqnarray}\label{+}
b_i^+\chi_i(S)=\frac{-x_i+\sqrt{x_i^2+4|\chi_i(S)|^2}}{2}a_i^+,~~~a_i^+\overline{\chi_i(S)}=\frac{x_i+\sqrt{x_i^2+4|\chi_i(S)|^2}}{2}b_i^+,
\end{eqnarray}
and
\begin{eqnarray}\label{-}
b_i^-\chi_i(S)=\frac{-x_i-\sqrt{x_i^2+4|\chi_i(S)|^2}}{2}a_i^-,~~~a_i^-\overline{\chi_i(S)}=\frac{x_i-\sqrt{x_i^2+4|\chi_i(S)|^2}}{2}b_i^-
\end{eqnarray}
where $x_i=\chi_i(R)-\chi_i(L)\in\Bbb R$. If $\chi_i(S)=0$, then we may assume that $\lambda_i^+=\chi_i(R)$, $\lambda_i^-=\chi_i(L)$, $(a_i^+,b_i^+)=(1,0)$ and $(a_i^-,b_i^-)=(0,1)$.  Let $c_i^+=\frac{|a_i^+|^2}{|a_i^+|^2+|b_i^+|^2}$, $c_i^-=\frac{|a_i^-|^2}{|a_i^-|^2+|b_i^-|^2}$, $d_i^+=\frac{|b_i^+|^2}{|a_i^+|^2+|b_i^+|^2}$, $d_i^-=\frac{|b_i^-|^2}{|a_i^-|^2+|b_i^-|^2}$, $e_i^+=\frac{a_i^+\overline{b_i^+}}{|a_i^+|^2+|b_i^+|^2}$ and $e_i^-=\frac{a_i^-\overline{b_i^-}}{|a_i^-|^2+|b_i^-|^2}$.  If $\chi_i(S)\neq 0$, then
\begin{eqnarray*}
&&c_i^+=\frac{|x_i+\sqrt{x_i^2+4|\chi_i(S)|^2}|^2}{|x_i+\sqrt{x_i^2+4|\chi_i(S)|^2}|^2+4|\chi_i(S)|^2},\\
&&c_i^-=\frac{|x_i-\sqrt{x_i^2+4|\chi_i(S)|^2}|^2}{|x_i-\sqrt{x_i^2+4|\chi_i(S)|^2}|^2+4|\chi_i(S)|^2},\\
&&d_i^+=\frac{4|\chi_i(S)|^2}{|x_i+\sqrt{x_i^2+4|\chi_i(S)|^2}|^2+4|\chi_i(S)|^2},\\
&&d_i^-=\frac{4|\chi_i(S)|^2}{|x_i-\sqrt{x_i^2+4|\chi_i(S)|^2}|^2+4|\chi_i(S)|^2},\\
&&e_i^+=2\chi_i(S)\frac{x_i+\sqrt{x_i^2+4|\chi_i(S)|^2}}{|x_i+\sqrt{x_i^2+4|\chi_i(S)|^2}|^2+4|\chi_i(S)|^2},\\
&&e_i^-=-2\chi_i(S)\frac{x_i+\sqrt{x_i^2+4|\chi_i(S)|^2}}{|x_i+\sqrt{x_i^2+4|\chi_i(S)|^2}|^2+4|\chi_i(S)|^2},.
\end{eqnarray*}

Furthermore, by \cite[Lemma 11]{AT},
\[v_i^+=\frac{1}{\sqrt{n(|a_i^+|^2+|b_i^+|^2)}}(a_i^+\chi_i(g_1^{-1}),\ldots,a_i^+\chi_i(g_n^{-1}),b_i^+\chi_i(g_1^{-1}),\ldots,b_i^+\chi_i(g_n^{-1}))^T\]
and
\[v_i^-=\frac{1}{\sqrt{n(|a_i^-|^2+|b_i^-|^2)}}(a_i^-\chi_i(g_1^{-1}),\ldots,a_i^-\chi_i(g_n^{-1}),b_i^-\chi_i(g_1^{-1}),\ldots,b_i^-\chi_i(g_n^{-1}))^T\]
, where $x^T$ is the transpose of vector $x$, are eigenvectors of $\Gamma$ corresponding to $\lambda_i^+$ and $\lambda_i^-$, respectively, and $\{v_i^+,v_i^-\mid 1\leq i\leq n\}$ is an orthonormal basis of $\Bbb C^{2n}$.  Let $E_i^+=v_i^+.(v_i^+)^*$ and $E_i^-=v_i^-.(v_i^-)^*$. Then we have
\begin{eqnarray*}
E_i^+=\frac{1}{n(|a_i^+|^2+|b_i^+|^2)}
\begin{bmatrix}
  |a_i^+|^2B&a_i^+\overline{b_i^+} B  \\
  b_i^+\overline{a_i^+}B &|b_i^+|^2B
\end{bmatrix},~~~~E_i^-=\frac{1}{n(|a_i^-|^2+|b_i^-|^2)}
\begin{bmatrix}
  |a_i^-|^2B&a_i^-\overline{b_i^-} B  \\
  b_i^-\overline{a_i^-}B &|b_i^-|^2B
\end{bmatrix},
\end{eqnarray*}
where $B=[\chi_i(g_r^{-1}g_s)]_{1\leq r,s\leq n}$. We keep these notations throughout the paper.

\section{PST on semi-Cayley graphs over abelian groups}\label{sec:pst}

Let us first give some relations between PST on Cayley graphs and semi-Cayley graphs.
\begin{lemma} Let $\Gamma=\SC(G,R,L,\emptyset)$. Then $\Gamma$ has a PST from $(g,r)$ to $(h,s)$ if and only if $r=s=0$ and $\Cay(G,R)$ has a PST from $g$ to $h$, or $r=s=1$ and $\Cay(G,L)$ has a PST from $g$ to $h$.
\end{lemma}
\proof Let $A$, $B$ and $C$ be the adjacency matrices of $\Gamma$, $\Cay(G,R)$ and $\Cay(G,L)$, respectively. Then
$
A=\begin{bmatrix}
B&0\\
0&C
\end{bmatrix}
$, which implies that $H_A(t)=\begin{bmatrix}
H_B(t)&0\\
0&H_C(t)
\end{bmatrix}$. Now the result is straightforward. 
\eqd

The proof of the following lemma is similar to \cite[Theorem 1]{Canul}. We give it for completeness.
\begin{lemma} Let $X=\Cay(G,R)$ and $\Gamma=\SC(G,R,R,S)$, where $G$ is an abelian group. Let $A$, $B$ and $C$ be the adjacency matrices of $\Gamma$, $X$ and $\Cay(G,S)$ (may be directed, or having loops), respectively. Then
\[H_A(t)=\begin{bmatrix}
 H_B(t)D_1&C D_2\\
 C^TD_2&H_B(t)D_1
 \end{bmatrix},\]
 where $D_1=\sum_{k=0}^\infty\frac{(-{\bf i}t)^{2k}}{(2k)!}(CC^T)^k=\cos(t\sqrt{CC^T})$ and $D_2=\sum_{k=0}^\infty\frac{(-{\bf i}t)^{2k+1}}{(2k+1)!}(CC^T)^k=-{\bf i}(CC^T)^{\frac{-1}{2}}\sin(t\sqrt{CC^T})$ where the second equality for $D_2$ holds if 
 $\det(C)\neq 0$. Furthermore, if $X$ has a PST from $g$ to $h$ at time $t$, then
 \begin{itemize}
 \item[(1)] $\Gamma$ has a PST from $(g,0)$ to $(h,0)$ and from $(g,1)$ to $(h,1)$ at time $t$, whenever $t\in\Bbb Z\pi$ and $CC^T=I$.
 \item[(2)] $\Gamma$ has a PST from $(g,0)$ to $(h,1)$ at time $t$, whenever $t\in(2\Bbb Z+1)\frac{\pi}{2}$ and $S=\{1\}$.
 \end{itemize}
\end{lemma}
\proof Suppose that  $G=\{g_1,g_2,\ldots,g_n\}$ and consider the ordering 
\[(g_1,0),\ldots,(g_n,0),(g_1,1),\ldots,(g_n,1)\]
 for the vertices of $\Gamma$. Then $A=I_2\otimes B+(P\otimes C+Q\otimes C^{T})$, where $I_2$ is the $2\times 2$ identity matrix, $P=\begin{bmatrix}
0&1\\
0&0
\end{bmatrix}$ and $Q=\begin{bmatrix}
0&0\\
1&0
\end{bmatrix}$. Since $G$ is abelian, $CC^{T}=C^{T}C$, $BC=CB$ and $BC^{T}=C^TB$. Also we have $P^2=Q^2=0_2$ and $PQ+QP=I_2$, where $0_2$ is the $2\times 2$ all $0$ matrix. Let $D=P\otimes C+Q\otimes C^{T}$. Then 
\[(I_2\otimes B)D=P\otimes (BC)+Q\otimes (BC^T)=P\otimes(CB)+Q\otimes (C^TB)=D(I_2\otimes B),\]
 and so $\exp(-{\bf i}tA)=\exp(-{\bf i}t(I_2\otimes B))\exp(-{\bf i}tD)=(I_2\otimes\exp(-{\bf i}tB))\exp(-{\bf i}tD)$. On the other hand, 
 \[D^{2k}=I_2\otimes (CC^T)^k,~~D^{2k+1}=P\otimes C(CC^T)^{k}+Q\otimes C^T(CC^T)^k,~~\forall ~k\geq 0.\]
 This implies that 
 \[H_A(t)=\exp(-{\bf i}tA)=\begin{bmatrix}
 \exp(-{\bf i}tB)D_1&C D_2\\
 C^TD_2&\exp(-{\bf i}tB)D_1
 \end{bmatrix}=\begin{bmatrix}
 H_B(t)D_1&C D_2\\
 C^TD_2&H_B(t)D_1
 \end{bmatrix},\]
 where $D_1=\sum_{k=0}^\infty\frac{(-{\bf i}t)^{2k}}{(2k)!}(CC^T)^k=\cos(t\sqrt{CC^T})$ and $D_2=\sum_{k=0}^\infty\frac{(-{\bf i}t)^{2k+1}}{(2k+1)!}(CC^T)^k=-{\bf i}(CC^T)^{\frac{-1}{2}}\sin(t\sqrt{CC^T})$ where the second equality for $D_2$ holds if 
 $\det(C)\neq 0$.
\eqd

In the rest of paper, we study the PST of semi-Cayley graphs in general using spectra of these graphs. Let $\Gamma=\SC(G,R,L,S)$ be a semi-Cayley graph over an abelian group $G$, $A$ be its adjacency matrix and 
$H(t)=\exp(-{\bf i}tA)$. By the notations of the previous section and the Spectral Decomposition Theorem, we have
\[A=\lambda_1^+E_1^++\lambda_1^-E_1^-+\ldots+\lambda_n^+E_n^++\lambda_n^-E_n^-,\]
and therefore we have the decomposition of the transfer matrix as follows:
\[H(t)=\exp(-{\bf i}\lambda_1^+t)E_1^++\exp(-{\bf i}\lambda_1^-t)E_1^-+\ldots+\exp(-{\bf i}\lambda_n^+t)E_n^++\exp(-{\bf i}\lambda_n^-t)E_n^-.\]

Let $u=(g,r)$ and $v=(h,s)$ be two vertices of $\Gamma$, where $r,s\in\{0,1\}$. Then
\begin{eqnarray}\label{main}
H_{u,v}(t)=\left\{
\begin{array}{ll}\frac{1}{n}\sum_{i=1}^n\Big(c_i^+\exp(-{\bf i}\lambda_i^+t)+c_i^-\exp(-{\bf i}\lambda_i^-t)\Big)\chi_i(g^{-1}h)&r=s=0\\
\frac{1}{n}\sum_{i=1}^n\Big(d_i^+\exp(-{\bf i}\lambda_i^+t)+d_i^-\exp(-{\bf i}\lambda_i^-t)\Big)\chi_i(g^{-1}h)&r=s=1\\
\frac{1}{n}\sum_{i=1}^n\Big(e_i^+\exp(-{\bf i}\lambda_i^+t)+e_i^-\exp(-{\bf i}\lambda_i^-t)\Big)\chi_i(g^{-1}h)&r=0,s=1\\
\frac{1}{n}\sum_{i=1}^n\Big(\overline{e_i^+}\exp(-{\bf i}\lambda_i^+t)+\overline{e_i^-}\exp(-{\bf i}\lambda_i^-t)\Big)\chi_i(g^{-1}h)&r=1,s=0
\end{array}
\right.
\end{eqnarray}

\begin{lemma}\label{base1} Keeping the above notations, we have
\begin{itemize}
\item[(i)] if $\chi_i(S)=0$ then $c_i^+=d_i^-=1$, $c_i^-=d_i^+=e_i^+=e_i^-=0$ and if $\chi_i(S)\neq 0$ then $e_i^++e_i^-=0$, and 
$0<c_i^+,c_i^-,d_i^+,d_i^-<1$. Furthermore, for each $i$, we have $c_i^++c_i^-=1$, $0< d_i^++d_i^-\leq 1$,
$|e_i^+|=|e_i^-|\leq \frac{1}{2}$ and so $|e_i^+|+|e_i^-|\leq 1$.
\item[(ii)] Then for all $i$ with $\chi_i(S)\neq 0$, we have $e_i^+\neq e_i^-$ and
 \begin{eqnarray*}
 c_i^+=c_i^-\Leftrightarrow d_i^+=d_i^-\Leftrightarrow x_i=0
 \end{eqnarray*}
\item[(iii)] $|H_{u,v}(t)|\leq 1$ for all $u,v$. 
\end{itemize}
\end{lemma}
\proof By a tedious computation, one can proof $(i)$ and $(ii)$.  Since $G$ is abelian, for each $g\in G$ and $\chi\in\Irr(G)$ we have $|\chi(g)|=1$. So the last part follows from the first part and equality (\ref{main}).
\eqd

Let $z_1,\ldots,z_k\in\Bbb C$. Then, it is easy to see that $|z_1+\ldots+z_k|=|z_1|+\ldots+|z_k|$ if and only if $z_j$'s have
the same argument. We use this fact in the following lemma.
\begin{lemma}\label{2sided}
Let $G$ be a finite abelian group of order $n$, $\Irr(G)=\{\chi_1,\ldots,\chi_n\}$, where $\chi_1$ is the trivial character of $G$, $\Gamma=\SC(G,R,L,S)$ and $X=\{i\mid \chi_i(S)=0, 1\leq i\leq n\}$. Let $u=(g,r)$, $v=(h,s)$ be two vertices and $a=g^{-1}h$.  Then $\Gamma$ has a PST from $u$ to $v$ at time $t$ if and only if one of the following holds
\begin{itemize}
\item[(1)] $r=s=0$ and 
\begin{eqnarray*}
\chi_i(a)&=&\exp(-{\bf i}(\lambda_1^+-\lambda_i^+)t)~~~\forall i=1,\ldots,n,\\
\chi_j(a)&=&\exp(-{\bf i}(\lambda_1^+-\lambda_j^-)t)~~~\forall j\notin X,
\end{eqnarray*}
\item[(2)] $r=s=1$ and
\begin{eqnarray*}
\chi_i(a)&=&\exp(-{\bf i}(\lambda_1^--\lambda_i^-)t)~~~\forall i=1,\ldots,n,\\
\chi_j(a)&=&\exp(-{\bf i}(\lambda_1^--\lambda_j^+)t)~~~\forall j\notin X,
\end{eqnarray*}
\item[(3)] $r=0, s=1$, $X=\emptyset$ and for all $j=1,\ldots,n$
\begin{eqnarray*}
\chi_j(a)&=&\frac{|\chi_j(S)|}{\chi_j(S)}\exp(-{\bf i}(\lambda_1^+-\lambda_j^+)t),\\
\exp(-{\bf i}\lambda_j^+t)&=&-\exp(-{\bf i}\lambda_j^-t),
\end{eqnarray*}
and $R=L$,
\item[(4)] $r=1, s=0$, $X=\emptyset$ and for all $j=1,\ldots,n$
\begin{eqnarray*}
\chi_j(a)&=&\frac{|\chi_j(S)|}{\overline{\chi_j(S)}}\exp(-{\bf i}(\lambda_1^+-\lambda_j^+)t), \\
\exp(-{\bf i}\lambda_j^+)t&=&-\exp(-{\bf i}\lambda_j^-t)
\end{eqnarray*}
and $R=L$.
\end{itemize}
\end{lemma}
\proof  By Lemma \ref{base1}, $|H_{u,v}(t)|\leq 1$ for all $u,v$. Let $\Gamma$ has a PST from $u$ to $v$.
Then $|H_{u,v}(t)|=1$. Now the argument  of all terms of the summation given in (\ref{main}) are the same number. First suppose that $r=s=0$. Then
\begin{eqnarray*}
\exp(-{\bf i}\lambda_i^+t)\chi_i(a)=\exp(-{\bf i}\lambda_1^+t)=\exp(-{\bf i}\lambda_j^-t)\chi_j(a), ~\forall i=1,\ldots,n,~\forall j\notin X,
\end{eqnarray*}

and therefore 
\begin{eqnarray*}
\chi_i(a)&=&\exp(-{\bf i}(\lambda_1^+-\lambda_i^+)t)~~~i=1,\ldots,n,\\
\chi_j(a)&=&\exp(-{\bf i}(\lambda_1^+-\lambda_j^-)t)~~~j\notin X.
\end{eqnarray*}
  
Similarly, one can proof the result for
the case $r=s=1$.  

Now let $r=0$ and $s=1$. Then $H_{u,v}(t)=\frac{1}{n}\sum_{j\notin X}\Big(e_j^+\exp(-{\bf i}\lambda_j^+t)+e_j^-\exp(-{\bf i}\lambda_j^-t)\Big)\chi_j(g^{-1}h)$. Since $|e_j^+|+|e_j^-|\leq 1$, we have $1=|H_{u,v}(t)|\leq \frac{1}{n}(n-|X|)$, which means that $X=\emptyset$, i.e for each $j=1,\ldots,n$, $\chi_j(S)\neq 0$. Furthermore, $e_j^-=-e_j^+$, and so
\[H_{u,v}(t)=\frac{1}{n}\sum_{j=1}^n e_j^+\Big(\exp(-{\bf i}\lambda_j^+t)-\exp(-{\bf i}\lambda_j^-t)\Big)\chi_j(g^{-1}h).\]

For each $j$, let $\theta_j^+$ be the argument of $e_j^+$. Then $\exp({\bf i}\theta_j^+)=\frac{\chi_j(S)}{|\chi_j(S)|}$. Again, since $|H_{u,v}(t)|=1$, the argument  of all terms of $H_{u,v}(t)$ are the same number. Hence
\[\arg(e_j^+\exp(-{\bf i}\lambda_j^+t)\chi_j(g^{-1}h))=\arg(-e_j^+\exp(-{\bf i}\lambda_j^-t)\chi_j(g^{-1}h))=\arg(\exp(-{\bf i}\lambda_1^+t)),\]
which implies that
 \[e_j^+\Big(\exp(-{\bf i}\lambda_j^+t)-\exp(-{\bf i}\lambda_j^-t)\Big)\chi_j(g^{-1}h)=0,~~~\chi_j(g^{-1}h)=\frac{|\chi_j(S)|}{\chi_j(S)}\exp(-{\bf i}(\lambda_1^+-\lambda_j^+)).\]
We have $\chi_j(g^{-1}h)\neq 0$ and also $e_j^+\neq 0$  since $\chi_j(S)\neq 0$. Thus the first equality implies that 
$\exp(-{\bf i}\lambda_j^+t)=-\exp(-{\bf i}\lambda_j^-t)$. Hence 
\[H_{u,v}(t)=\frac{1}{n}\sum_{j=1}^n 2e_j^+\exp(-{\bf i}\lambda_j^+t)\chi_j(g^{-1}h)=\frac{1}{n}\sum_{j=1}^n 2e_j^+\frac{|\chi_j(S)|}{\chi_j(S)}.\]
Since $|2e_j^+\frac{|\chi_j(S)|}{\chi_j(S)}|\leq 1$ and $|H_{u,v}(t)|=1$, for all $j=1,\ldots,n$, we have $e_j^+=\frac{|\chi_j(S)|}{2\chi_j(S)}$, which implies that $x_j=0$. Thus we have proved that $\chi_j(R)=\chi_j(L)$ for all $j=1,\ldots,n$. Now, by the column orthogonality relations of characters \cite[p. 161]{JL}, for all $g\in G$ we have
\begin{eqnarray*}
\sum_{i=1}^n\chi_i(g^{-1})\chi_i(R)=\sum_{i=1}^n\sum_{x\in R}\chi_i(g^{-1})\chi_i(x)=\sum_{x\in R}\sum_{i=1}^n\chi_i(g^{-1})\chi_i(x)=|G|\sum_{x\in R}\delta_{xg},
\end{eqnarray*}
and similarly $\sum_{i=1}^n\chi_i(g^{-1})\chi_i(L)=|G|\sum_{y\in L}\delta_{yg}$, where $\delta$ is the characteristic function. Therefore, $\sum_{x\in R}\delta_{xg}=\sum_{y\in L}\delta_{yg}$ for all $g\in G$. This implies that $R=L$.

 To prove the last part, it is enough to note that we have to consider the complex conjugates of $e_j^+$ and $e_j^-$ is the summation. The converse direction is straightforward. This completes the proof.
\eqd

\begin{remark}
Since $1\notin X$ if and only if $S\neq\emptyset$, by Lemma \ref{2sided} if $S\neq\emptyset$ and $\Gamma$ has a PST at time $t$ between two
pairs $u=(g,r), v=(h,s)$ and $u'=(g',r'),v'=(h',s')$, then $r=s$ if and only if $r'=s'$. 
\end{remark}

\begin{corollary}
Let $\Gamma=\SC(G,R,L,S)$ be a semi-Cayley graph over an abelian group $G$ and $\Gamma$ has a PST between two vertices $(g,r)$ and $(h,s)$, with $r\neq s$. Then $R=L$ and $\Gamma$ is a Cayley graph over a group isomorphic to $G\rtimes\Bbb Z_2$.
\end{corollary}
\proof By Lemma \ref{2sided}, we have $R=L$. Since $G$ is abelian, the map $\tau:G\rightarrow G$, by the rule $\tau(g)=g^{-1}$
is an automorphism of $G$. Then \cite[Lemma 3.2]{ZF} implies that $\Gamma$ is a Cayley graph over a group isomorphic to $G\rtimes\Bbb Z_2$.
\eqd

The following result, characterizes Cayley graphs over finite groups with an abelian groups of index $2$ having PST.
\begin{corollary}
Let $\Gamma=\Cay(G,T)$ be an undirected Cayley graph, where $G$ has an abelian subgroup $H$ of index $2$, $G=H\cup xH$, $T=T_1\cup xT_2$, where $T_1,T_2\subset H$ (if $T_2=\emptyset$ then we put $xT_2=\emptyset$), $\Irr(H)=\{\chi_1,\ldots,\chi_n\}$ and $X=\{i\mid \chi_i(T_2)=0\}$. Then eigenvalues of $\Gamma$ are
\[\lambda_i^+=\frac{\chi_i(T_1)+\chi_i(xT_1x^{-1})+\sqrt{(\chi_i(T_1)-\chi_i(xT_1x^{-1}))^2+4|\chi_i(x^2T_2)|}}{2},~~i=1,\ldots,n\]
and 
\[\lambda_i^-=\frac{\chi_i(T_1)+\chi_i(xT_1x^{-1})-\sqrt{(\chi_i(T_1)-\chi_i(xT_1x^{-1}))^2+4|\chi_i(x^2T_2)|}}{2},~~i=1,\ldots,n.\]
Furthermore, $\Gamma$ has a PST
between two vertices $g_1$ and $g_2$ at time $t$ if and only if one of the following holds
\begin{itemize}
\item[(1)] $g_1,g_2\in H$, and 
\begin{eqnarray*}
\chi_i(a)&=&\exp(-{\bf i}(\lambda_1^+-\lambda_i^+)t)~~~\forall i=1,\ldots,n,\\
\chi_j(a)&=&\exp(-{\bf i}(\lambda_1^+-\lambda_j^-)t)~~~\forall j\notin X,
\end{eqnarray*}
where $a=g_1^{-1}g_2$.
\item[(2)] $x^{-1}g_1,x^{-1}g_2\in H$, and 
\begin{eqnarray*}
\chi_i(a)&=&\exp(-{\bf i}(\lambda_1^--\lambda_i^-)t)~~~\forall i=1,\ldots,n,\\
\chi_j(a)&=&\exp(-{\bf i}(\lambda_1^--\lambda_j^+)t)~~~\forall j\notin X,
\end{eqnarray*}
where $a=g_1^{-1}g_2$.
\item[(3)] $g_1\in H$, $x^{-1}g_2\in H$, and for all $j=1,\ldots,n$, $\chi_j(x^2T_2)\neq 0$ and 
\begin{eqnarray*}
\chi_j(a)&=&\frac{|\chi_j(x^2T_2)|}{\chi_j(x^2T_2)}\exp(-{\bf i}(\lambda_1^+-\lambda_j^+)t), \\
\exp(-{\bf i}(\lambda_j^+)t)&=&-\exp(-{\bf i}\lambda_j^-t),
\end{eqnarray*}
and $T_1x=xT_1$, where $a=(xg_1)^{-1}g_2$.
\item[(4)] $x^{-1}g_1\in H$, $g_2\in H$ and for all $j=1,\ldots,n$, $\chi_j(x^2T_2)\neq 0$ and
\begin{eqnarray*}
\chi_j(a)&=&\frac{|\chi_j(x^2T_2)|}{\overline{\chi_j(x^2T_2)}}\exp(-{\bf i}(\lambda_1^+-\lambda_j^+)t), \\
\exp(-{\bf i}(\lambda_j^+)t)&=&-\exp(-{\bf i}\lambda_j^-t),
\end{eqnarray*}
and $T_1x=xT_1$, where $a=(x^{-1}g_1)^{-1}g_2$.
\end{itemize}
\end{corollary}
\proof By \cite[Lemma 8]{AT}, we may assume that $G=H\cup xH$, $\Gamma=\SC(H,R,L,S)$ and we can identify $(h,0)$ and $(h,1)$ with $h$ and $xh$, respectively, where $R=H\cap T=T_1$, $L=H\cap xTx^{-1}=xT_1x^{-1}$ and $S=H\cap xT=x^2T_2$. Since $T=T^{-1}$ and $T_1\cap xT_2=\emptyset$, we see that $T_1^{-1}=T_1$ and $T_2^{-1}=T_2$. Now the result is a direct consequence of Lemma \ref{2sided}.
\eqd

As the following result shows, if a semi-Cayley graph $\SC(G,R,L,S)$, with $S=S^{-1}$, over an abelian group $G$ of order $n$ has a PST between distinct vertices 
then $n$ is even. In particular, if a Cayley graph over an abelian group of order $2n$ , or over the dihedral group of order $2n$ has a PST then $n$ is even. This proves \cite[Theorem 3.5]{TFC} and  also \cite[Theorem 8]{CCL} follows from this and the last part of Theorem \ref{R=L}.  
\begin{corollary}\label{odd} Let $\Gamma=\SC(G,R,L,S)$, $\Gamma$ has a PST between two distinct vertices $(g,r)$ and $(h,s)$ and put $a=g^{-1}h$. If $r=s$ then $a$ has order $2$, and if $r\neq s$ then $a$ has order $2$ if and only if $S=S^{-1}$. In particular, if $G$ has odd order then $\Gamma$ has no PST between two distinct vertices $(g,r)$ and $(h,r)$, and moreover if $S=S^{-1}$ then $G$ has no PST between any two distinct vertices.
\end{corollary}
\proof
Let $a=g^{-1}h$ has order $m$. Then 
$m\neq 1$, since $a\neq 1$. Then for each $i$, $\chi_i(a)$ is an $m$-th root of unity, since $G$ is abelian. Thus we may assume that for each $i$, $\chi_i(a)=\exp(-{\bf i}\frac{2\pi}{m}k_{a,i})$, where $k_{a,i}\in\Bbb Z$.
Moreover, we may assume that $(m,k_{a,i})=1$. Let $\chi$ be the complex conjugate of $\chi_i$.
Then $\chi=\chi_j$ for some $j$. Then $\lambda_j^+=\lambda_i^+$ since $R$ and $L$ are inverse-closed. Furthermore, $k_{a,j}=-k_{a,i}$. Then, Lemma \ref{2sided} implies that if $r=s=0$ then
 $\frac{k_{a,i}}{m}-(\lambda_1^+-\lambda_i^+)T\in\Bbb Z$ and $-\frac{k_{a,i}}{m}-(\lambda_1^+-\lambda_i^+)T\in\Bbb Z$
for all $i=1,\ldots,n$, 
where $T=\frac{t}{2\pi}$. So for each $i$, we have $2\frac{k_{a,i}}{m}\in\Bbb Z$. Therefore $m=2$ since $(k_{a,i},m)=1$. Similarly, one can prove the result for the case $r=s=1$. 

Now suppose that $r\neq s$. If $S=S^{-1}$ then $\overline{\chi(S)}=\chi(S)$ for all $\chi\in\Irr(G)$, which means that $\chi(S)$ is a real number for all $\chi\in\Irr(G)$. 
Now, by the same notations and applying the same proof of the first part, Lemma \ref{2sided} implies that for each $j=1,\ldots,n$, we have $2\frac{k_{a,j}}{m}\in\Bbb Z$, where $(k_{a,j},m)=1$. This means that $m=2$. Conversely, suppose that $m=2$. Then
for all $j$, $\chi_j(a)=1$ or $-1$. If $r=0$ and $s=1$ then Lemma \ref{2sided}, $\frac{\chi_j(S)}{|\chi_j(S)|}=\pm \exp(-{\bf i}(\lambda_1^+-\lambda_j^+)t)$ and also by considering the conjugate of $\chi_j$, we have $\frac{\overline{\chi_j(S)}}{|\chi_j(S)|}=\pm \exp(-{\bf i}(\lambda_1^+-\lambda_j^+)t)$, which implies that $\chi_j(S)=\overline{\chi_j(S)}=\chi_j(S^{-1})$ for all $j=1,\ldots,n$. Again, by the same argument using the column orthogonality relations given in the proof of Lemma \ref{2sided}, $S=S^{-1}$. This completes the proof.
\eqd

Now we are ready to give a characterization of $\SC(G,R,L,S)$ having PST, where $G$ is abelian and $R=L$.  For this purpose, we  need notation of the $2$-adic exponential valuation of rational numbers which is a mapping defined by
\begin{align*}
\begin{split}
&\psi_2:\mathbb{Q}\rightarrow\mathbb{Z}\cup\{\infty\},\\
&\psi_2(0)=\infty,\\
&\psi_2(2^l\frac{a}{b})=l,~~~~where~ a,b,l \in \mathbb{Z} ~and~ 2\not|ab.\\
\end{split}
\end{align*}

\noindent We assume that $\infty+\infty=\infty+l=\infty$ and $\infty> l$ for any $l \in \mathbb{Z}$. Then $\psi_2$ has the following properties. For $\beta, \beta' \in \mathbb{Q}$,
\begin{itemize}
  \item[(1)] $\psi_2(\beta\beta')=\psi_2(\beta)+\psi_2(\beta'),$
  \item[(2)] $\psi_2(\beta+\beta') \geq {\rm min}(\psi_2(\beta), \psi_2(\beta'))$ and the equality holds if   $\psi_2(\beta)\ne \psi_2(\beta').$
\end{itemize}
\begin{theorem}\label{R=L} Let $G$ be a finite abelian group of order $n$, $\Irr(G)=\{\chi_1,\ldots,\chi_n\}$, where $\chi_1$
 is the trivial character of $G$, and $\Gamma=\SC(G,R,R,S)$. Let $u=(g,r)$ and $v=(g,r)$ are vertices of $\Gamma$ and $a=g^{-1}h$. If $r=s$ then $\Gamma$ has a PST between two distinct vertices  
$u$ and $v$ if and only if the following  conditions hold
\begin{itemize}
\item[(i)] $a$ has order $2$,
 \item[(ii)] $\Gamma$ is integral, and for each $i$, $\chi_i(R)$ and $|\chi_i(S)|$ are integers.
 \item[(iii)] $\nu_2(|R|+|S|-\lambda_i^+)=\nu_2(|R|+|S|-\lambda_i^-)$ is the same integer, say  $k$, for all $i$ that $\chi_i(g_x^{-1}g_y)=-1$ and for all $i$ with $\chi_i(g_x^{-1}g_y)=1$, $\nu_2(|R|+|S|-\lambda_i^+)>k$ and $\nu_2(|R|+|S|-\lambda_i^-)>k$.
\end{itemize}
Also $\Gamma$ is periodic if and only if  $\Gamma$ is integral. Furthermore, the minimum
period of the vertices is $\frac{2\pi}{M}$, where $M=\gcd\{\lambda-\lambda_1^+\mid \lambda\in\Spec(\Gamma)\setminus\{\lambda_1^+\}\}$.
\end{theorem}
\proof 
  Let $r=s$. Since $R=L$, we have $\lambda_i^+=\chi_i(R)+|\chi_i(S)|$ and $\lambda_i^-=\chi_i(R)-|\chi_i(S)|$. Let $X=\{i\mid \chi_i(S)=0\}$. If $j\in X$ then $\lambda_j^+=\lambda_j^-$. Then by Lemma \ref{2sided}, $\Gamma$ has a PST between vertices $u$ and $v$ if and only if for all $i$ we have $\exp(-{\bf i}\lambda_i^+t)\chi_i(a)=\exp(-{\bf i}\lambda_i^-t)\chi_i(a)$.
Let $T=\frac{t}{2\pi}$. Then we have
$(\lambda_i^+ -\lambda_i^-)T=2|\chi_i(S)|T\in\Bbb Z$ for all $i$.
If we put $i=1$ in the above relation, then $2|S|T\in\Bbb Z$ which means that $T\in\Bbb Q$. Thus $|\chi_i(S)|\in\Bbb Q$ for all $i$. Therefore $|\chi_i(S)|\in\Bbb Z$ for all $i$, since they are eigenvalues of $\SC(G,\emptyset,\emptyset,S)$.

By the same argument of the proof of Corollary \ref{odd}, for all $i$ we have
\[\frac{k_{a,i}}{m}-(\lambda_1^+-\lambda_i^+)T, \frac{k_{a,i}}{m}-(\lambda_1^+-\lambda_i^-)T\in\Bbb Z, \]
where $a$ has order $m$ and $k_{a,i}$ is an integer.
Since $T, \frac{k_{a,i}}{m}\in\Bbb Q$, we have $\lambda_1^+-\lambda_i^+\in\Bbb Q$ and $\lambda_1^+-\lambda_i^-\in\Bbb Q$. So $\lambda_i^+,\lambda_i^-\in\Bbb Q$.
Therefore $\Gamma$ is an integral graph, since $\lambda_i^+$ and $\lambda_i^-$ are algebraic integers. Hence for each $i$, $\chi_i(R)\in\Bbb Q$.  Thus $\chi_i(R)\in\Bbb Z$ for all $i$, since $\chi_i(R)$, $i=1,\ldots,n$ are the eigenvalues of  Cayley graph $\Cay(G,R)$. 

Again by the proof of Corollary \ref{odd}, if $u$ and $v$ are distinct then $m=2$ and so $\chi_i(a)=1$ or $-1$ for all $i$.
First assume that $\chi_i(a)=-1$. Then
$(\lambda_1^+-\lambda_i^+)T, (\lambda_1^+-\lambda_i^-)T\in\Bbb Z+\frac{1}{2}$. So
$\nu_2(\lambda_1^+-\lambda_i^+)+\nu_2(T)=\nu_2((\lambda_1^+-\lambda_i^+)T)=-1$ and $\nu_2(\lambda_1^+-\lambda_i^-)+\nu_2(T)=\nu_2((\lambda_1^+-\lambda_i^-)T)=-1$. Thus $\nu_2(\lambda_1^+-\lambda_i^+)=\nu_2(\lambda_1^+-\lambda_i^-)=-1-\nu_2(T)=k$. Now assume that
$\chi_i(a)=1$. Then
 $(\lambda_1^+-\lambda_i^+)T, (\lambda_1^+-\lambda_i^-)T\in\Bbb Z$ and by similar argument  $\nu(|R|+|S|-\lambda_i^+)=k+l_1$ and $\nu(|R|+|S|-\lambda_i^-)=k+l_2$ for some $l_1,l_2\geq 1$. This proves the first part. The second part is a direct consequence
 of the proof of the first part.
\eqd

\begin{theorem}\label{added}
Let $\Gamma=\SC(G,R,L,S)$, and $u=(g,r), v=(h,s)$, where $r\neq s$ be two distinct vertices of $\Gamma$.  Then $\Gamma$ has 
a PST between $u$ and $v$ at time $t$ if and only if the following conditions hold
\begin{itemize}
\item[(i)] $R=L$,
\item[(ii)] $\chi_j(S)\neq 0$ for each $j$,
\item[(iii)] if $r=0, s=1$ then $\chi_j(g^{-1}h)=\frac{|\chi_j(S)|}{\chi_j(S)} \exp(-{\bf i}(\lambda_1^+-\lambda_j^+)t)$,
and if $r=1,s=0$ then $\chi_j(g^{-1}h)=\frac{|\chi_j(S)|}{\overline{\chi_j(S)}} \exp(-{\bf i}(\lambda_1^+-\lambda_j^+)t)$, where in both cases $\exp(-{\bf i}(\lambda_1^+-\lambda_j^+)t)=\pm 1$.
\item[(iv)] $\chi_j(R), |\chi_j(S)|\in\Bbb Z$ for each $j$, in particular $\Gamma$ is integral,
\item[(v)] $v_2(|S|)=v_2(|\chi_j(S)|)$ for all $j$.
\end{itemize}
\end{theorem}
\proof Suppose that $\Gamma$ has a PST between $u$ and $v$ at time $t$ and $T=\frac{t}{2\pi}$. Then $R=L$, for each $j$ we have $\chi_j(S)\neq 0$ and
$\exp(-{\bf i}\lambda_j^+t)=-\exp(-{\bf i}\lambda_j^-t)$ by Lemma \ref{2sided}.  The last equality implies that
$(\lambda_j^+-\lambda_j^-)T-\frac{1}{2}\in\Bbb Z$. Hence $2|\chi_j(S)|T=(\lambda_j^+-\lambda_j^-)T\in\Bbb Q$, for all $j$. Putting $j=1$, we get $T\in\Bbb Q$ and so 
$\lambda_j^+-\lambda_j^-\in\Bbb Q$ for all $j$. Hence $|\chi_j(S)|\in\Bbb Q$ for all $j$. Since the eigenvalues of 
$\SC(G,\emptyset,\emptyset,S)$ are $\pm |\chi_j(S)|$, we conclude that $|\chi_j(S)|\in\Bbb Z$ for all $j$. Furthermore, $2|\chi_j(S)|T-\frac{1}{2}\in\Bbb Z$ imples that $4|\chi_j(S)|T$ is an odd number and so $0=v_2(4|\chi_j(S)|T)=2+v_2(|\chi_j(S)|)+v_2(T)$ for all $j$. Therefore $v_2(|S|)=v_2(|\chi_j(S)|)=-2-v_2(T)$ for all $j$.

On the other hand, again by Lemma \ref{2sided}, if $r=0, s=1$ or $r=1,s=0$ then for each $j$ we have $\chi_j(a)=\frac{|\chi_j(S)|}{\chi_j(S)}\exp(-{\bf i}(\lambda_1^+-\lambda_j^+)t)$ or $\chi_j(a)=\frac{|\chi_j(S)|}{\overline{\chi_j(S)}}\exp(-{\bf i}(\lambda_1^+-\lambda_j^+)t)$, respectively. Since the complex conjugate of any irreducible character of $G$ is an irreducible character of $G$ and the corresponding
eigenvalues of an irreducible character and its complex conjugate are the same, since $R$ is inverse-closed, we conclude that in both cases  $\exp(-{\bf i}(\lambda_1^+-\lambda_j^+)t)$ is a real number, which means that it is $1$ or $-1$. Hence $(\lambda_1^+-\lambda_j^+)T\in\frac{1}{2}\Bbb Z\subset\Bbb Q$. Since $T\in\Bbb Q$ and $\lambda_1^+=|R|+|S|\in\Bbb Z$, we have $\lambda_j^+\in Q$. This means that $\Gamma$ is integral and $\chi_j(R)\in\Bbb Z$ for all $j$. The converse direction is clear by Lemma \ref{2sided}.
\eqd

\begin{remark} By Theorem \ref{R=L} and Corollary \ref{odd}, one can exactly determine the existence of a PST between vertices of an 
undirected Cayley graph
 $\Cay(G,T)$, whenever $G$ has
an abelian sugroup of index $2$ and $H\cap T=H\cap xTx^{-1}$, where $\{1,x\}$ is a left transversal to $H$ in $G$. In
particular, since if $G$ is an abelian group of even order, a dihedral group, a dicyclic group or a semidihedral group, then it satisfies the
 above conditions, one
can reprove the main results of \cite{ASGh}, \cite{CCL}, \cite{CF}, and \cite{TFC}. For example, 
let $\Gamma=\Cay(G,T)$, where $G$ is an abelian group of order $2n$. Then there exists a subgroup $H\leq G$, which is abelian and of index $2$. Let $G=H\cup xH$ and $T=T_1\cup xT_2$, where $T_1,T_2\subseteq H$. Then $xT_1x^{-1}=T_1$ and so $\Gamma=\SC(H,T_1,T_1,x^2T_2)$. Let $\Irr(H)=\{\chi_1=1,\ldots,\chi_n\}$. Then by Theorem \ref{R=L},  if $1\neq g_1^{-1}g_2\in H$  then $\Gamma$ has a PST between  $g_1$ and $g_2$ 
if and only if the following  conditions hold
\begin{itemize}
\item[(i)] the order of $g_1^{-1}g_2$ is two,
 \item[(ii)] $\Gamma$ is integral, and for each $i$, $\chi_i(T_1)$ and $|\chi_i(T_2)|$ are integers,
 \item[(iii)] $\nu_2(|T|-\lambda_i^+)=\nu_2(|T|-\lambda_i^-)$ is the same integer, say  $k$, for all $i$ that $\chi_i(g_1^{-1}g_2)=-1$ and for all $i$ with $\chi_i(g_1^{-1}g_2)=1$, $\nu_2(|T|-\lambda_i^+)>k$ and $\nu_2(|R|+|S|-\lambda_i^-)>k$,
\end{itemize}
where 
$\lambda_i^+=\chi_i(T_1)+|\chi_i(T_2)|$ and $\lambda_i^-=\chi_i(T_1)-|\chi_i(T_2)|$, $i=1,\ldots,n$.

Also, if $g_1^{-1}g_2\notin H$, then $\Gamma$ has a PST between  $g_1$ and $g_2$ at time $t$
if and only if the following  conditions hold
\begin{itemize}
\item[(i)] $\chi_j(T_2)\neq 0$ for each $j$,
\item[(ii)] if $g_1\in H$ and $g_2\notin H$ then $\chi_j(g_{1}^{-1}g_2)=\frac{|\chi_j(T_2)|}{\chi_j(x^2T_2)} \exp(-{\bf i}(\lambda_1^+-\lambda_j^+)t)$,
and if $g_1\notin H$ and $g_2\in H$ then $\chi_j(g_1^{-1}g_2)=\frac{|\chi_j(T_2)|}{\overline{\chi_j(x^2T_2)}} \exp(-{\bf i}(\lambda_1^+-\lambda_j^+)t)$, where in the both cases $\exp(-{\bf i}(\lambda_1^+-\lambda_j^+)t)=\pm 1$.
\item[(ii)] $\chi_j(T_1), |\chi_j(T_2)|\in\Bbb Z$ for each $j$, in particular $\Gamma$ is integral,
\item[(vi)] $v_2(|T_2|)=v_2(|\chi_j(T)|)$ for all $j$.
\end{itemize}

Also $\Gamma$ is periodic if and only if  $\Gamma$ is integral. Furthermore, the minimum
period of the vertices is $\frac{2\pi}{M}$, where $M=\gcd\{\lambda-\lambda_1^+\mid \lambda\in\Spec(\Gamma)\setminus\{\lambda_1^+\}\}$. This proves that \cite[Theorem 2.4]{TFC}, the main theorem of \cite{TFC}, is a direct consequence of Theorems \ref{R=L} and \ref{added}.
\end{remark}

Let us recall generalized dihedral and generalized dicyclic groups as follows: Given an abelian group $A$, the generalized dihedral group $Dih(A,x)=\langle A,x\rangle$ is a group generated by $A$ and an
 element $x$ such that $x\notin A$, $x^2=1$ and $xax=a^{-1}$ for all $a\in A$. It is easy to see that $A$ is a normal subgroup
 of $Dih(A)$ of index $2$. Clearly the dihedral group $D_{2n}=\langle a,b\rangle\mid a^n=b^2=(ab)^2=1\rangle$ of order $2n$,
  is $Dih(\langle a\rangle,b)$. Now let $A$ be an abelian group of even order and of exponent greater than $2$, and let $y$ be an involution of $A$. The
  generalized dicyclic group $Dic(A,y,x)$ is the group $\langle A,x\mid x^2=y,x^{-1}ax=a^{-1},~ \forall a\in A\rangle$. When
  $A$ is cyclic, $Dic(A,y,x)$ is called a dicyclic (or generalized quaternion) group.
  
The following result improves \cite[Theorem 8 and Theorem 11]{CCL} and \cite[Theorem 4.2]{ASGh}.
\begin{corollary}\label{gen}
Let $A$ be an abelian group, $\Irr(A)=\{1=\chi_1,\ldots,\chi_n\}$, $G=Dih(A,x)$ or $G=Dic(A,y,x)$, and $\Gamma=\Cay(G,T)$, where $T=T_1\cup xT_2$ for some $T_1,T_2\subseteq A$. if $1\neq g_1^{-1}g_2\in A$  then $\Gamma$ has a PST between  $g_1$ and $g_2$ 
if and only if the following  conditions hold
\begin{itemize}
\item[(i)] the order of $g_1^{-1}g_2$ is two,
 \item[(ii)] $\Gamma$ is integral, and for each $i$, $\chi_i(T_1)$ and $|\chi_i(T_2)|$ are integers,
 \item[(iii)] $\nu_2(|T|-\lambda_i^+)=\nu_2(|T|-\lambda_i^-)$ is the same integer, say  $k$, for all $i$ that $\chi_i(g_1^{-1}g_2)=-1$ and for all $i$ with $\chi_i(g_1^{-1}g_2)=1$, $\nu_2(|T|-\lambda_i^+)>k$ and $\nu_2(|R|+|S|-\lambda_i^-)>k$,
\end{itemize}
where 
$\lambda_i^+=\chi_i(T_1)+|\chi_i(T_2)|$ and $\lambda_i^-=\chi_i(T_1)-|\chi_i(T_2)|$, $i=1,\ldots,n$.

Also, if $g_1^{-1}g_2\notin A$, then $\Gamma$ has a PST between  $g_1$ and $g_2$ at time $t$
if and only if the following  conditions hold
\begin{itemize}
\item[(i)] $\chi_j(T_2)\neq 0$ for each $j$,
\item[(ii)] if $g_1\in A$ and $g_2\notin A$ then $\chi_j(g_{1}^{-1}g_2)=\frac{|\chi_j(T_2)|}{\chi_j(T_2)} \exp(-{\bf i}(\lambda_1^+-\lambda_j^+)t)$,
and if $g_1\notin A$ and $g_2\in A$ then $\chi_j(g_1^{-1}g_2)=\frac{|\chi_j(T_2)|}{\overline{\chi_j(T_2)}} \exp(-{\bf i}(\lambda_1^+-\lambda_j^+)t)$, where in the both cases $\exp(-{\bf i}(\lambda_1^+-\lambda_j^+)t)=\pm 1$.
\item[(ii)] $\chi_j(T_1), |\chi_j(T_2)|\in\Bbb Z$ for each $j$, in particular $\Gamma$ is integral,
\item[(vi)] $v_2(|T_2|)=v_2(|\chi_j(T)|)$ for all $j$.
\end{itemize}

Also $\Gamma$ is periodic if and only if  $\Gamma$ is integral. Furthermore, the minimum
period of the vertices is $\frac{2\pi}{M}$, where $M=\gcd\{\lambda-\lambda_1^+\mid \lambda\in\Spec(\Gamma)\setminus\{\lambda_1^+\}\}$. 
\end{corollary}
\proof Since $T_1=A\cap T=A\cap xTx^{-1}$ and $T_2=A\cap xT$, we have $\Gamma=\SC(A,T_1,T_1,T_2)$. Now the result follows from Theorems \ref{R=L} and  \ref{added} immediately.
\eqd

Since $\SC(G,R,L,G)$ is the join graph of $\Cay(G,R)$ and $\Cay(G,L)$, the following result gives a characterization of the join of two Cayley graphs over the same abelian group. 
\begin{corollary}\label{S=G}
Let $\Gamma=\SC(G,R,L,G)$, where $G$ is an abelian group of order $n$, and $\Irr(G)=\{\chi_1,\ldots,\chi_n\}$. Then $\Gamma$ has a PST between vertices $u=(g,r)$ and $v=(h,s)$ at time $t$ if and only if one of the following holds:
\begin{itemize}
\item[(1)] $g=h, r=s=0$, $(\lambda_1^+-\chi_i(R))t\in 2\pi\Bbb Z$ for all $i$ and 
$t\sqrt{(|R|-|L|)^2+4n^2}\in\pi\Bbb Z$.
\item[(2)] $g=h, r=s=1$, $(\lambda_1^--\chi_i(L))t\in 2\pi\Bbb Z$ for all $i$ and 
$t\sqrt{(|R|-|L|)^2+4n^2}\in\pi\Bbb Z$.
\item[(3)] $g\neq h, r=s=0$, $\sum_{i=1}^n\exp({\bf i}\lambda_i^+t)=0$ and $t\sqrt{(|R|-|L|)^2+4n^2}\in\pi\Bbb Z$.
\item[(4)] $g\neq h, r=s=1$, $\sum_{i=1}^n\exp({\bf i}\lambda_i^-t)=0$ and 
$t\sqrt{(|R|-|L|)^2+4n^2}\in\pi\Bbb Z$,
\end{itemize}
where $\lambda_1^+=\frac{|R|+|L|+|\sqrt{(|R|-|L|)^2+4n^2}}{2}$, $\lambda_1^-=\frac{|R|+|L|-|\sqrt{(|R|-|L|)^2+4n^2}}{2}$
, $\lambda_i^+=\chi_i(R)$ and $\lambda_i^-=\chi_i(L)$, $i=2,\ldots,n$ are eigenavlues of $\Gamma$.
\end{corollary}
\proof We have $X=\{i\mid \chi_i(S)=0\}=\{2,\ldots,n\}$ since $\chi(G)=0$ for all $1\neq \chi\in\Irr(G)$. This menas that for all $i\neq 1$, we may assume that $\lambda_i^+=\chi_i(R)$ and $\lambda_i^-=\chi_i(L)$. Then, by Lemma \ref{2sided} and Theorem \ref{added}, $\Gamma$ has a PST from $u=(g,r)$ to $v=(h,s)$ if and only if $r=s=0$, $\chi_i(a)=\exp(-{\bf i}(\lambda_1^+-\lambda_i^+)t)$ for all $i$ and $t\sqrt{(|R|-|L|)^2+4n^2}\in\pi\Bbb Z$ or $r=s=1$,  $\chi_i(a)=\exp(-{\bf i}(\lambda_1^--\lambda_i^-)t)$ for all $i$ and $t\sqrt{(|R|-|L|)^2+4n^2}\in\pi\Bbb Z$ where $a=g^{-1}h$. Putting $a=1$, we get 
the first and second parts. Now let $a\neq 1$. Then by the column orthogonality relations \cite[Theorem 16.4]{JL}, we have $\sum_{i=1}^n\chi_i(a)=0$, which proves the last two parts.
\eqd

Now we give some examples as an application of our results as follows:

\begin{example} Let $n\geq 3$ and $\Gamma$ be the $n$-sunlet graph, a graph that its vertices is obtained by attaching $n$ pendant edges to the $n$-cycle. Then $\Gamma=\SC(G,R,L,S)$, where $G=\langle a\rangle\cong\Bbb Z_n$, $R=\{a,a^{-1}\}$, $L=\emptyset$ and $S=\{1\}$. If $n$ is odd, then $\Gamma$ has no PST by Corollary \ref{odd}. If $n=2m$ is even, by considering the character $\chi_m:G\rightarrow\Bbb C$ by the rule $\chi_m(a^k)=(-1)^k$, Lemma \ref{2sided} implies that $\Gamma$ has no PST and is not periodic.
\end{example}

\begin{example} Let $n\geq 3$ and $\Gamma$ be the $n$-gonal $n$-cone, the join of an $n$-cycle with the empty graph with $n$ vertices. 
Then $\Gamma=\SC(G, R,L,S)$, where $G=\langle a\rangle\cong\Bbb Z_n$, $R=\{a,a^{-1}\}$, $L=\emptyset$ and $S=G$. Then, 
by Corllary \ref{S=G}, $\lambda_1^+=1+\sqrt{1+n^2}$, $\lambda_1^-=1-\sqrt{1+n^2}$, and for all $i\neq 1$, we have $\lambda_i^+=2Re(\chi_i(a))=2\cos(\frac{2\pi i}{n})$ and $\lambda_i^-=0$. If $n$ is odd, then $\Gamma$ has no PST by Corollary \ref{odd}.
Let $n=2m$. Consider the irreducible character $\chi_m:G\rightarrow\Bbb C$ with the rule $\chi_m(a^k)=(-1)^k$ and let
$\Gamma$ has a PST from $(g,r)$ to $(h,s)$. Then Corollary \ref{S=G} implies that $g=h, r=s=0$ or $g=h,r=s=1$ occurs if $\sqrt{1+n^2}\in\Bbb Q$, which is impossible. If $g\neq h$ and $r=s=1$ then we must have $\exp({\bf i}t)=1-n-(-{\bf i})^k$ for some integer $k$, which is again impossible. Also, if $g\neq h$ and $r=s=0$ then we must have $({\bf i})^k\exp({\bf i}t)+\sum_{i=1}^{n-1}\exp(2{\bf i}\cos(\frac{2\pi i}{n}))=0$, for some integer $k$.
\end{example}

\begin{example}
Let $\Gamma=\Cay(G,T)$, where $G=Dih(A,x)$ or $G=Dic(A,y,x)$ and $T=xA$. Then, by the notations of Corollary \ref{gen}, we have  
$T_1=\emptyset$, $T_2=A$, $\lambda_1^+=|A|$, $\lambda_1^-=-|A|$, and for all $i\neq 1$, $\lambda_i^+=\lambda_i^-=0$. Then $\Gamma$ is periodic with minimum period $\frac{2\pi}{|A|}$ and has no PST between distinct vertices, by Corollary \ref{gen}.
\end{example}


\begin{example}
Let $\Gamma=\Cay(G,T)$, where $G=Dih(A,x)$ or $G=Dic(A,y,x)$ and $T=xA\cup \{a\in A\mid o(a)= 2\}$. Then, by the notations of Corollary \ref{gen}, $T_1=\{a\in A\mid o(a)=2\}$, $T_2=A$, $\lambda_1^+=|T_1|+|A|=|T|$, $\lambda_1^-=|T_1|-|A|$, and for each
$i\neq 1$, $\lambda_i^+=\lambda_i^-=\sum_{a\in A, o(a)=2}\chi_i(a)$. Again, $\Gamma$ is integral and so it is periodic with
minimum period $\frac{2\pi}{M}$, where $M=\gcd(|T_1|-|A|-|T|,\sum_{a\in A, o(a)=2}\chi_2(a)-|T|, \ldots, \sum_{a\in A, o(a)=2}\chi_n(a)-|T|)$, and $n=|A|$. Also by a tedious computation, we see that $\Gamma$ has no PST.
\end{example}

\section{Concluding remarks}
The perfect state transfer (PST) in quantum walks, where the underlying graphs are semi-Cayley graphs over abelian groups were studied.  We derive the spectral decomposition of any semi-Cayley graph over an abelian group $G$ in terms of the irreducible characters  of $G$. Since every Cayley graph over a group $G$ having an abelian subgroup $H$ is a semi-Cayley graph over $H$,
our method for investigating of PST can be applied to any Cayley graph over a finite group having an abelian subgroup of index $2$. One can apply our method to the quasi-abelian semi-Cayley graphs.

\textbf{Acknowledgements}
The author would like to thank Tao Feng (faculty member of Department of Mathematics, Beijing Jiaotong University). In fact one day after sumbitting the first version of my paper to Arxiv.org, Tago Feng 
informed me that there is a mistake in my paper and the equality $|e_i^+|+|e_i^-|\leq \frac{1}{2}$ is not correct. This caused me to
revise the article and  fix the problem in Theorem \ref{added}. Also note that in the previous version on arXiv.org (2202.03062v2)
I forgot to add the acknowledgement.


\begin{thebibliography}{Dillo 150}
\bibitem{AT} M. Arezoomand and B. Taeri, On the characteristic polynomial of $n$-Cayley digraphs, Electron. J. Combin., 20(3)
(2013), \# P57.
\bibitem{AT3} M. Arezoomand and B. Taeri, A classification of finite groups with integral bi-Cayley graphs, Trans. Combin.,4(4) (2015) 55-61.
\bibitem{ASGh} M. Arezoomand, F. Shafiei and Modjtaba Ghorbani, Perfect state transfer on Cayley graphs over the dicyclic group, Linear Algebra Appl., 639 (2022) 116-134.
\bibitem{Canul} R.J. Angeles-Canul, R. Norton, M. Opperman, C. Paribello, M. Russell, C. Tamon, Perfect state transfer, integral circulants and join of graphs, Quantum
Computation and Information 10 (2010) 325–342. arXiv:0907.2148.
\bibitem{CCL} X. Cao, B. Chen and S. Ling, Perfect state transfer on Cayley graphs over the dihedral groups: the non-normal case, Electron. J. Combin., 27 (2020) no.2 P2.28. 
\bibitem{CF} X. Cao, and K. Feng, Perfect state transfer on Cayley graphs over dihedral groups, Linear and Multilinear Algebra, https://doi.org/10.1080/03081087.2019.1599805.
\bibitem{CGG} G. Coutinho\ et al., Perfect state transfer on distance-regular graphs and association schemes, Linear Algebra Appl.  478 (2015) 108--130.
\bibitem{C}  G.  Coutinho,  Quantum state transfer in graphs [PhD dissertation]. University of Waterloo, 2014.
\bibitem{FG} E. Farhi, S. Gutmann, Quantum computation and decision trees, Phys. Rev. A 58 (3) (1998) 915–928.
 
\bibitem{G1} C. Godsil, Periodic graphs, Electron. J. Combin.  18 (2011) no.~1 Paper 23, 15 pp.

\bibitem{G2} C. Godsil, State transfer on graphs, Discrete Math.  312 (2012) no.~1 129--147.

\bibitem{G3} C. Godsil, When can perfect state transfer occur?, Electron. J. Linear Algebra 23 (2012) 877--890.

\bibitem{HL} J. Huang and S. Li, Integral and distance integral Cayley graphs over generalized dihedral groups, J. Algebraic Combin., 
https://doi.org/10.1007/s10801-020-00948-1.

\bibitem{JL} G. James and M. Liebeck, Representations and characters of groups, Cambridge University Press, Second Edition, 2001.
\bibitem{J}    N. Johnston, et al., Perfect quantum state transfer using Hadamard diagonalizable graphs, Linear Algebra Appl.  531 (2017) 375--398.

\bibitem{LCWW} G. Luo, X. Cao, D. Wang and X. Wu, Perfect quantum state transfer on Cayleygraphs over semi-dihedral groups, Linear and Multilinear Algebra, https://doi.org/10.1080/03081087.2021.1954585.

\bibitem{SS}   M. \v{S}tefa\v{n}\'{a}k, S. Skoup\'{y}, Perfect state transfer by means of discrete-time quantum walk on complete bipartite graphs, Quantum Inf. Process.  16 (2017) no.~3 Paper No. 72, 14 pp.


\bibitem{TFC} Y.Y. Tan, K. Feng, X. Cao, Perfect state transfer on abelian Cayley graphs, Linear Algebra Appl., (2019), https://doi.org/10.1016/j.laa.2018.11.011 

\bibitem{matrix} F. Zhang, Matrix Theory, Basic Results and Techniques, Second Edition, Springer, 2011. 
\bibitem{ZF} J. X Zhou and Y. Q. Feng, cubic bi-Cayley graphs over abelian groups, European J. Combin. 36 (2014) 679-793.






























\end{thebibliography}
\end{document}